\documentclass[12pt]{amsart}
\usepackage{ifthen,verbatim}
\usepackage{graphicx}
\usepackage{mathrsfs}
\usepackage{amsmath,amssymb,amsthm}
\usepackage{tikz-cd}
\usepackage{color}
\usepackage{caption}
\usepackage{subcaption}
\numberwithin{equation}{section}
\nonstopmode
\theoremstyle{plain}

\newtheorem{thm}[equation]{Theorem}

\newtheorem{lem}[equation]{Lemma}

\theoremstyle{definition}
\newtheorem{rem}[equation]{Remark}

\newcommand{\RomanNumeralCaps}[1]
    {\MakeUppercase{\romannumeral #1}}

\newenvironment{pf}[1][]{%
 \vskip 3mm
 \noindent
 \ifthenelse{\equal{#1}{}}%
  {{\slshape Proof. }}%
  {{\slshape #1.} }%
 }%
{\qed\bigskip}

\newcounter{alphabet}

\newcommand{\be}{\begin{equation}}
\newcommand{\ee}{\end{equation}}
\newcommand{\ba}{\begin{align}}
\newcommand{\ea}{\end{align}}
\newcommand{\ben}{\begin{align*}}
\newcommand{\een}{\end{align*}}

\newcommand{\ve}{{\varepsilon}}

\newcommand{\K}{{\mathbf K}}

\newcommand{\R}{{\mathbb R}}
\newcommand{\C}{{\mathbb C}}

\newcommand{\uhp}{{\mathbb H}}
\newcommand{\Z}{{\mathbb Z}}

\newcommand{\Q}{{\mathbb Q}}

\newcommand{\sphere}{{\widehat{\mathbb C}}}
\renewcommand{\Im}{{\operatorname{Im}\,}}
\newcommand{\Aut}{{\operatorname{Aut}}}
\renewcommand{\Re}{{\operatorname{Re}\,}}

\newcommand{\id}{{\operatorname{id}}}
\newcommand{\fix}{{\operatorname{Fix}}}

\renewcommand{\mod}{\,{\operatorname{mod}\,}}

\newcommand{\inv}{^{-1}}
\newcommand{\area}{{\,\operatorname{Area}}}

\newcommand{\PSL}{{\operatorname{PSL}}}
\newcommand{\SL}{{\operatorname{SL}}}

\newcommand{\aand}{{\quad\text{and}\quad}}
\newcommand{\Gauss}{{\null_2F_1}}

\newcounter{minutes}
\divide\time by 60
\newcounter{hours}
\multiply\time by 60
\addtocounter{minutes}{-\time}

\begin{document}
\bibliographystyle{amsplain}
\title
{Geometric deduction of the solutions to modular equations}

\def\thefootnote{}
\footnotetext{
\texttt{\tiny File:~\jobname .tex,
          printed: \number\year-\number\month-\number\day,
          \thehours.\ifnum\theminutes<10{0}\fi\theminutes}
}
\makeatletter\def\thefootnote{\@arabic\c@footnote}\makeatother

\author[Md. S. Alam]{Md. Shafiul Alam}
\address{Graduate School of Information Sciences,
Tohoku University, Aoba-ku, Sendai 980-8579, Japan $-$ and $-$ Department of Mathematics, University of Barishal, Barishal-8254, Bangladesh}
\email{shafiul@dc.tohoku.ac.jp, msalam@bu.ac.bd}
\author[T. Sugawa]{Toshiyuki Sugawa}
\address{Graduate School of Information Sciences,
Tohoku University, Aoba-ku, Sendai 980-8579, Japan}
\email{sugawa@math.is.tohoku.ac.jp}

\begin{abstract}
In his notebooks, Ramanujan presented  without proof many remarkable 
formulae for the solutions to generalized modular equations.
Much later, proofs of the formulae were provided by making use of
highly nontrivial identities for theta series and hypergeometric functions.
We offer a geometric approach to the proof of those formulae.
We emphasize that our proofs are geometric and independent of such identities.
\end{abstract}
\thanks{The first author was supported by Data Sciences Program 
\RomanNumeralCaps{2} of Graduate School of Information Sciences, Tohoku University. 
The authors were supported in part by JSPS KAKENHI Grant Number JP17H02847.
}
\maketitle

\section{Introduction}

For  given integers $p\ge 2,$  Ramanujan, 
an Indian mathematical genius, considered the equation
\begin{equation}\label{eq:gme}
\frac{\Gauss(a,1-a;1;1-\beta)}{\Gauss(a,1-a;1;\beta)}
=p\frac{\Gauss(a,1-a;1;1-\alpha)}{\Gauss(a,1-a;1;\alpha)},
\end{equation}
known as the generalized modular equation of degree $p$ and signature $1/a.$
He left many formulae
describing relations between $\alpha$ and $\beta$ in his unpublished
notebooks but he did not record any proof of those formulae (see \cite{Berndt: NotebookIII} and \cite{Ram: lnb}).
Here, $\Gauss(a,b;c;z)$ denotes the hypergeometric function whose definition will be given
in Section 2.
The case when $1/a=2$ corresponds to the classical modular equation.
Indeed, the complete elliptic integral of the first kind is described by
$$
\K(r)=\int_0^1 \frac{dx}{\sqrt{(1-x^2)(1-r^2x^2)}}
=\frac\pi 2\Gauss(\tfrac12,\tfrac12;1;r^2)
$$
and the function
$$
\mu(r)=\frac\pi 2\cdot\frac{\K(\sqrt{1-r^2})}{\K(r)}
=\frac\pi 2\cdot\frac{\Gauss(\tfrac12,\tfrac12;1;1-r^2)}%
{\Gauss(\tfrac12,\tfrac12;1;r^2)}
$$
is known to be the modulus of the Gr\"otzsch ring $\{z\in\C: |z|<1\}
\setminus [0,r]$ for $0<r<1$ and plays an important role
in the theory of plane quasiconformal mappings (see, for instance, \cite{LV:qc}
or \cite{AVV:conf}).
The modular equation \eqref{eq:gme} for $1/a=2$ now takes the form 
$\mu(s)=p\mu(r)$ with $\alpha=r^2$ and $\beta=s^2.$
When $p=2,$ the solution to this modular equation is given by 
\begin{equation}\label{eq:classical}
s
=\frac{1-\sqrt{1-r^2}}{1+\sqrt{1-r^2}}
=\left(\frac{1-\sqrt{1-r^2}}{r}\right)^2
\end{equation}
(see \cite[(2.4) in p.~60]{LV:qc} or \cite[(5.4)]{AVV:conf}).

Ramanujan mainly considered the case when the signature $1/a$ is 3, 4 or 6.
See \cite{Berndt: NotebookI}, \cite{Berndt: NotebookII}, \cite{Berndt: NotebookIII}, 
\cite{Berndt: NotebookV}, \cite{BB:AGM} and \cite{Ram: nb} for his work
in relation with the modular equations.
Berndt, Bhargava and Garvan \cite{BBG95} gave proofs for those formulae of Ramanujan
(see also \cite{AQVV00}).
Their proofs make use of highly nontrivial identities for Jacobi's
theta functions and hypergeometric functions in addition to a number of ingeneous ideas.
For example, they gave rigorous proofs for the following results.

\begin{thm}[$\text{\cite[Theorem 7.1]{BBG95}}$]\label{thm:23}
When $p=2$ and $q=3,$ the solutions $\alpha$ and $\beta$
to the equation \eqref{eq:gme} are related by
\begin{equation}\label{eq:me23}
(\alpha\beta)^{1/3}+\big\{(1-\alpha)(1-\beta)\big\}^{1/3}=1.
\end{equation}
\end{thm}

\begin{thm}[$\text{\cite[Lemma 7.4]{BBG95}}$]\label{thm:33}
When $p=3$ and $q=3,$ the solutions $\alpha$ and $\beta$
to the equation \eqref{eq:gme} are related by
\begin{equation}\label{eq:me33}
(1-\alpha)^{1/3}=\frac{1-\beta^{1/3}}{1+2\beta^{1/3}}.
\end{equation}
\end{thm}

We note that the above relations can be transformed to polynomial equations.
For instance, \eqref{eq:me23} may be transformed to
\begin{equation}\label{eq:alg}
(2\alpha-1)^3\beta^3-3\alpha(4\alpha^2-13\alpha+10)\beta^2
+3\alpha(2\alpha^2-10\alpha+9)\beta-\alpha^3=0.
\end{equation}
In particular, we observe that there are at most three values of $\beta$
satisfying the modular equation \eqref{eq:me23} for each $\alpha.$
We can say that $\alpha$ and $\beta$ satisfy a polynomial equation 
of degree 3 in this case.
It is rather surprising that $\alpha$ and $\beta$ are related algebraically, because
the hypergeometric function is transcendental for the corresponding parameters.
For instance, we do not have a complete answer to the question for which $p$ and $a$ the
solutions to the modular equation \eqref{eq:gme} is algebraic.
In this note, we propose a geometric approach to this problem.
In particular, the geometric observation suggests that it is more natural
to look at $q=1/(1-2a)$ rather than the signature $1/a.$
We will call $q$ the \emph{order} of the modular equation \eqref{eq:gme}.
For instance, the signatures $1/a=2,3,4,6$ correspond to $q=\infty, 3, 2, 3/2,$
respectively.
Though our approach does not cover the case $1/a=6,$ it may allow us
to approach other cases when $q$ are integers $>3.$
Let $G_q=\langle T, W\rangle$ be the Fuchsian group generated by $T\tau=\tau+\lambda$
and $W\tau=\tau/(1+\lambda \tau),$ where $\lambda=2\cos\frac{\pi}{2q}$
and let $M_p$ be the M\"obius transformation defined by $M_p\tau=p\tau.$
To speak about the solutions to \eqref{eq:gme}, we have to clarify the range
of the solutions.
Originally, the range of $\alpha$ and $\beta$ was the interval $[0,1].$
However, in our setting, it is natural to choose the quotient Riemann surface 
$X=G_q\backslash\uhp$ as the range of the solutions.
Here, we take $X=\sphere\setminus\{0,1\}$ for finite $q$ and
$X=\C\setminus\{0,1\}$ for $q=\infty.$
Then we have the following result.

\begin{thm}\label{thm:Fuchs}
Let $n,\,p$ and $q$ be integers with $n,\,p,\,q\ge 2$ $($possibly $q=\infty)$
and set $a=(q-1)/(2q).$
Then, the solutions $(\alpha, \beta)$ to the equation \eqref{eq:gme}
in $X=G_q\backslash\uhp$
satisfy the equation $P(\alpha,\beta)=0$ for an irreducible polynomial
$P(x,y)$ of degree $n$ if and only if
$G_q\cap G_q^{M_p}$ is a subgroup of $G_q$ of index $n,$
where $G_q^{M_p}=M_p\inv G_q M_p.$
\end{thm}

As we will see in the proof, the solutions $\alpha$ and $\beta$
are parametrized as $\alpha=\varphi(z)$ and $\beta=\psi(z)$
on the Riemann surface $Z=(G_q\cap G_q^{M_p})\backslash\uhp$ and they
satisfy the polynomial equation $P(\alpha,\beta)=0$ of order $n.$
More precisely, $P(x,y)$ has the forms
\begin{align*}
P(x,y) &=a_0(x)y^n+a_1(x)y^{n-1}+\cdots +a_{n-1}(x)y+a_n(x) \\
&=b_0(y)x^n+b_1(y)x^{n-1}+\cdots +b_{n-1}(y)x+b_n(y),
\end{align*}
where $a_j(x)\in \C[x], b_j(y)\in\C[y]~(j=0,1,\dots, n),$
such that $P(x,y)$ is irreducible as an element of $\C(x)[y]$ and $\C(y)[x],$
respectively.
(Recall that $\C[x]$ and $\C(x)$ stand for the $\C$-algebra of polynomials in $x$
and the field of rational functions of $x$ with
coefficients in $\C,$ respectively.)

As applications of it, we will show \eqref{eq:classical} and Theorems
\ref{thm:23} and \ref{thm:33} in Sections 4 and 5.
In Section 2, we recall basic facts about hypergeometric functions
and Hecke groups and its subgroup of index 2.
Section 3 is devoted to the statement of our main lemma and its proof.
Theorem \ref{thm:Fuchs} is also proved in Section 3 with the help
of the main lemma.
We emphasize that our approach does not need ingenious ideas as in
\cite{BBG95}. Without knowledge about Ramanujan's formulae, we can
deduce those formulae possibly with laborious computations only.
In this paper, we treat only the case when the Riemann surface $Z$ is planar.
When $Z$ is non-planar, it is technically difficult to find an explicit form
of the polynomial $P(x,y).$
Let $\hat Z$ be the compactification of $Z.$
Since the parametrizations $\varphi,\psi:\hat Z\to\sphere$ are
$n$-sheeted covering maps with critical values contained in $\{0,1,\infty\}$
if $|G_q:G_q\cap G_q^{M_p}|=n,$ Belyi's theorem implies that
the compact Riemann surface $\hat Z$ is an algebraic curve defined 
over $\overline{\Q}$ (see \cite{JG:dessin}).
Therefore, in principle, we could determine the surface $Z$ and maps
$\varphi, \psi$ by the combinatorial information about the coverings.
We hope to give further examples when $Z$ is non-planar in the future work.

\medskip
\textbf{Acknowledgement.}
The authors would like to thank Professor Hiroshige Shiga for
suggestions which simplified the proof of Lemma \ref{lem:N} in the
early version of the manuscript.
They also thank Professor Matti Vuorinen for many suggestions.

\section{Hypergeometric functions and Hecke groups}\label{sec:basic}

The hypergeometric function $\Gauss(a,b;c;z)$ is defined by
$$
\Gauss(a,b;c;z)=\sum_{n=0}^\infty \frac{(a)_n (b)_n}{(c)_n n!}z^n,
\quad |z|<1,
$$
for parameters $a,b,c\in\C$ with $c\ne 0,-1,-2,\dots.$
By Euler's integral representation formula, we know that
$\Gauss(a,b;c;z)$ analytically extends to the slit domain
$\C\setminus [1,+\infty).$
Let $0<a\le 1/2.$
Then it is classically known \cite{Carat:f2} that the function
$$
\tau = f_a(z)=i\cdot\frac{\Gauss(a,1-a;1;1-z)}{\Gauss(a,1-a;1;z)}
$$
maps the upper half-plane $\uhp=\{z\in\C: \Im z>0\}$
onto the domain $\Delta_a=\{\tau\in\uhp: 0<\Re \tau<\sin\pi a,
|2\tau\sin\pi a-1|>1\}$ and that $\Delta_a$ is the hyperbolic triangle
with vertices at $0, \infty$ and $ie^{-\pi a i}$ of
interior angles $0, 0$ and $(1-2a)\pi,$ respectively 
(see also \cite[Lemma 4.1]{ASVV10}).
Note also that $f_a$ extends homeomorphically to the boundary
and 
\be\label{eq:fa}
f_a(0)=\infty,\quad f_a(1)=0,\quad f_a(\infty)=ie^{-\pi a i}=e^{\pi i(1-2a)/2}.
\ee
Suppose now that $q=1/(1-2a)$ is an integer or $\infty.$
Let $\pi_q:\Delta_a\to \uhp$ be the inverse map $z=f_a\inv(\tau)$
of $\tau=f_a(z).$
Then, by repeated applications of the Schwarz reflection principle,
$\pi_q$ extends to a holomorphic map from $\uhp$ into $\sphere\setminus\{0,1\}.$
Here, we note that $\pi_q$ is locally $q$ to 1 at the point $ie^{-\pi a i}.$
By construction, the covering group $G_q=\{\gamma\in\Aut(\uhp):
\pi_q\circ \gamma=\pi_q\}$ is the triangle group of signature $(q,\infty,\infty)$
arising from the hyperbolic triangle $\Delta_a,~ q=1/(1-2a).$
Here, the group $\Aut(\uhp)$ of analytic automorphism of $\uhp$ is identified
with $\PSL(2,\R)=\SL(2,\R)/\{\pm I\}.$
Since $G_q$ is a subgroup of the Hecke group $H_{2q}$ of index 2,
we recall basics of Hecke groups.
See \cite{CS98} for details about the Hecke groups.
Moreover, we refer to \cite{Beardon:disc} and \cite{Katok:fg} as general
references for Fuchsian groups.

By its form, the function $f_a$ satisfies the relation $f_a(1-z)=-1/f_a(z).$
Since $\pi_q=f_a\inv$ on $\Delta_a,$ the following result follows.

\begin{lem}\label{lem:recip}
The covering map $\pi_q$ satisfies the functional equation
$$
\pi_q(\tau)+\pi_q(-1/\tau)=1,\quad \tau\in\uhp.
$$
\end{lem}


For an integer $k\ge 3,$ the Hecke group $H_k$ is defined as the discrete subgroup
of $\PSL(2,\R)$ generated by the two elements $\pm S$ and $\pm T_k,$ where
$$
S=\begin{pmatrix}0&-1\\1&0\end{pmatrix},\quad
T_k=\begin{pmatrix}1&\lambda_k\\0&1\end{pmatrix}
\aand
\lambda_k=2\cos\frac\pi k.
$$
See \cite{Beardon:disc} or \cite{CS98}.
Here and in what follows, we often identify a matrix 
$A=\begin{pmatrix}a&b\\ c&d\end{pmatrix}$
with the M\"obius transformation $\tau\mapsto A\tau=(a\tau+b)/(c\tau+d)$
so that $S, T_k$ are regarded as elements of $\Aut(\uhp)=\PSL(2,\R).$
Note that $H_3$ is the classical modular group $\PSL(2,\Z).$
We put
$$
U_k=T_kS=\begin{pmatrix}\lambda_k&-1\\1&0\end{pmatrix}.
$$
Then $U_k$ is an elliptic element of order $k$ with fixed point at $e^{\pi i/k}.$
Let $F_k$ be the set of points $\tau\in\uhp$ with $|\Re \tau|\le \cos\frac\pi{k}, |\tau|\ge 1.$
Then $F_k$ is a fundamental domain for $H_k$ with side pairing transformations $S, T_k.$
(Note that $S$ has a fixed point at $\tau=i.$)
In particular, we see that $H_k$ is a triangle group of signature $(2,k,\infty).$

\begin{figure}[h]
     \centering
     \begin{subfigure}[b]{0.45\textwidth}
         \centering
         \includegraphics[width=\textwidth]{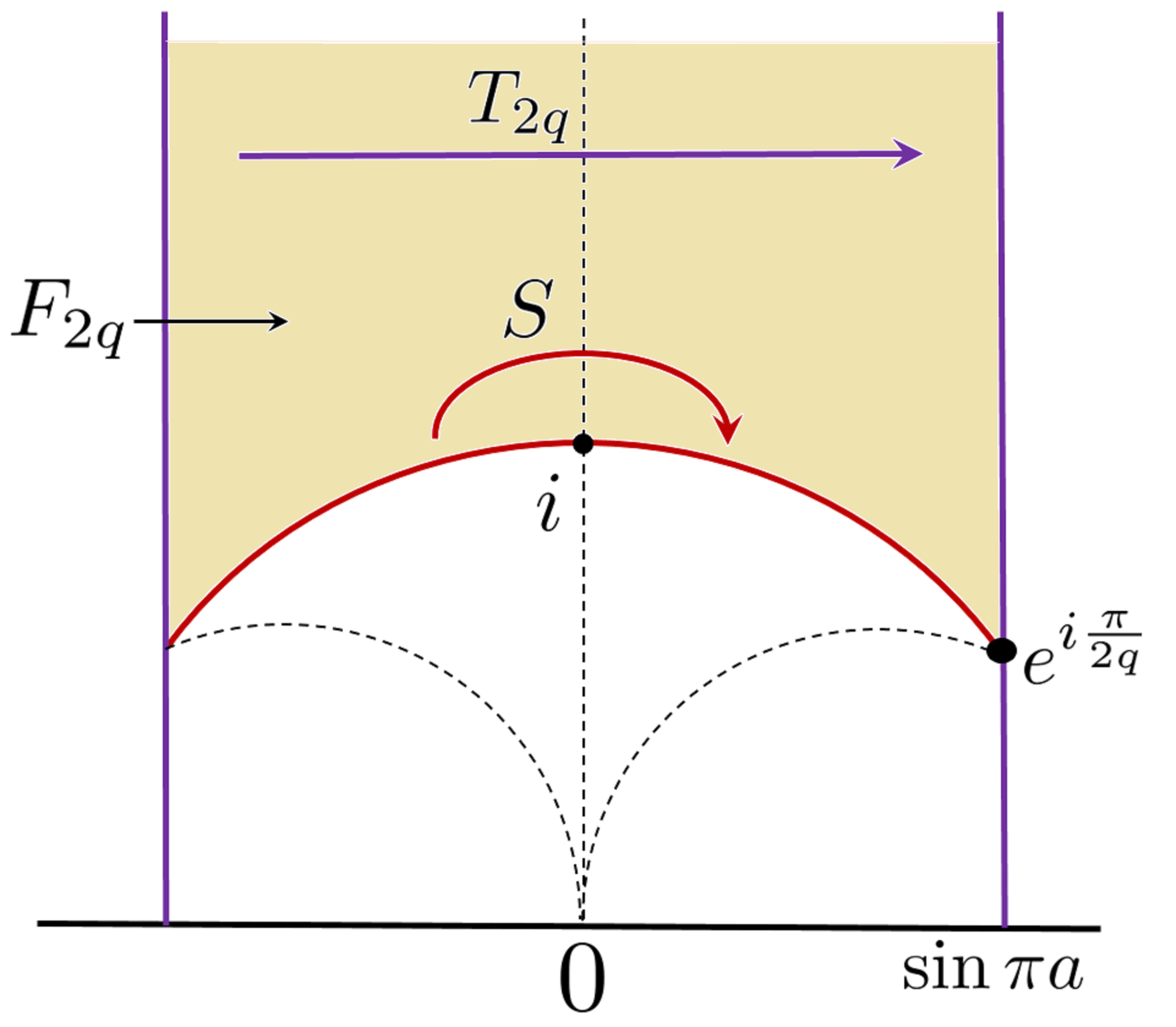}
     \end{subfigure}
     \hfill
     \begin{subfigure}[b]{0.45\textwidth}
         \centering
         \includegraphics[width=\textwidth]{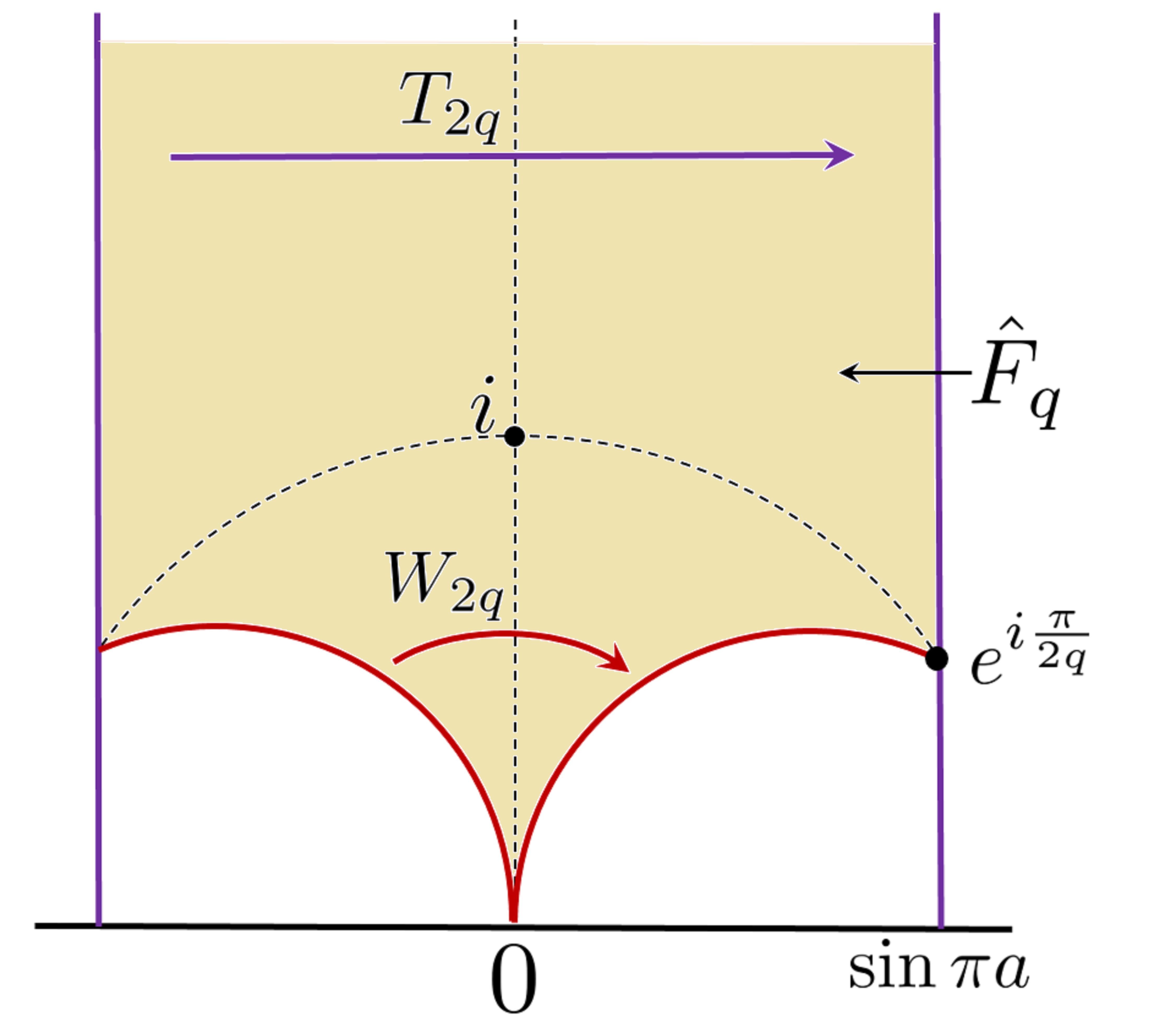}
     \end{subfigure}
           \caption{Fundamental domains of $H_{2q}=\langle T_{2q}, S\rangle$ and $G=\langle T_{2q}, W_{2q}\rangle$.}
        \label{fig:fq1}
\end{figure}

\begin{figure}[h]
     \centering
        \includegraphics[width=.6\textwidth]{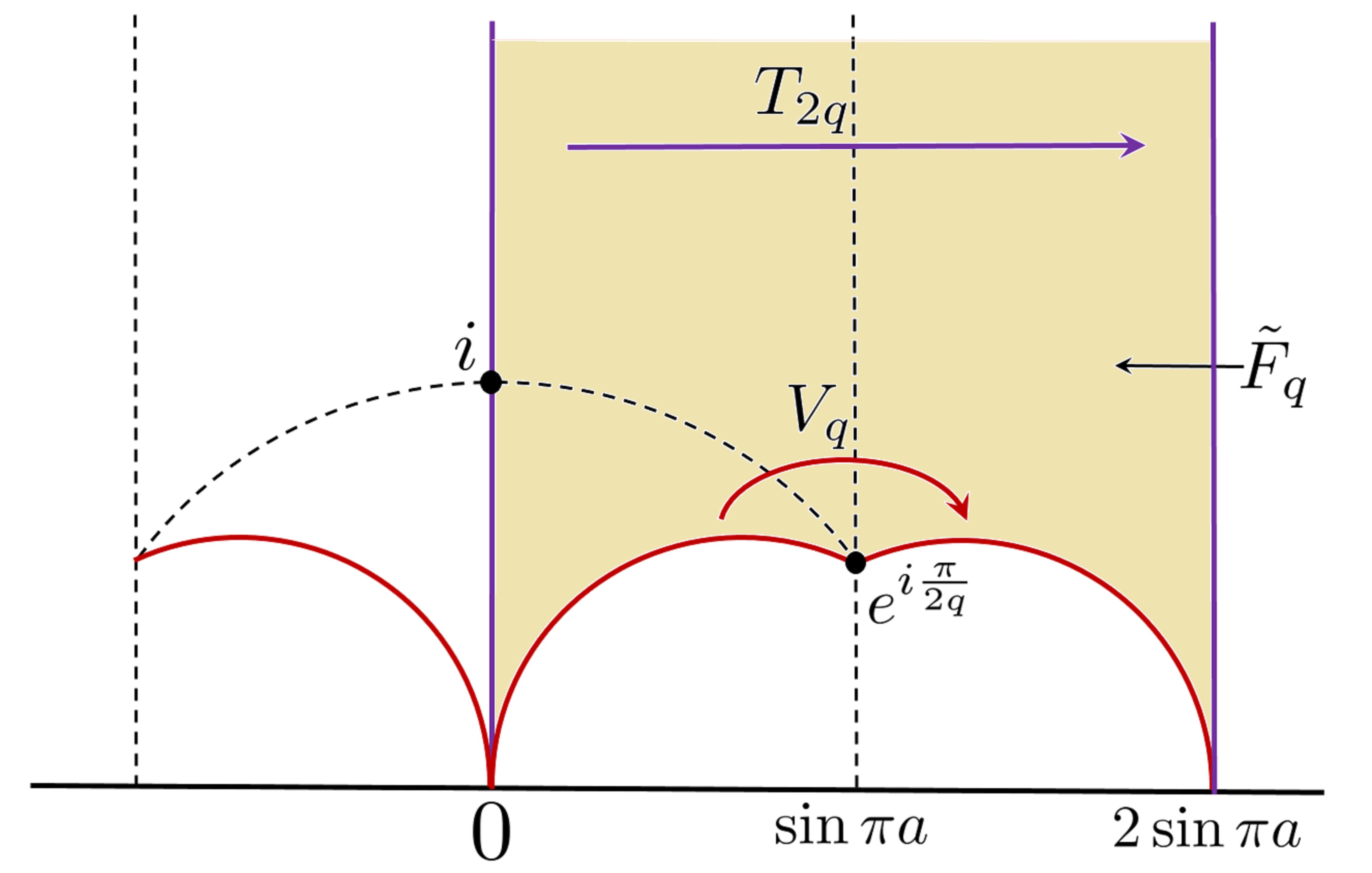}
         \caption{Fundamental domain $\tilde F_q$ of $G_q=\langle T_{2q}, V_q\rangle$.}
        \label{fig:fq2}
     \end{figure}

From now on, we restrict ourselves on the case when $k=2q$ is an even number
with $q\ge 2.$
Then $\hat F_q=F_{2q}\cup S(F_{2q})$ is a fundamental domain 
for the (normal) subgroup $G$ of $H_{2q}$ of index 2 generated by 
$$
T_{2q}=\begin{pmatrix}1&\lambda_{2q}\\0&1\end{pmatrix}
\aand
W_{2q}=S\inv T_{2q}\inv S=
\begin{pmatrix}1&0\\ \lambda_{2q}&1\end{pmatrix}.
$$
To adapt with our aim, we modify the fundamental domain as follows.
For $a=(q-1)/(2q),$ let $\tilde F_q=\Delta_a \cup \Delta_a'$, 
where $\Delta_a'$ is the reflection of $\Delta_a$ across the line $\Re \tau=\sin{\pi a}$. 
Then,  $\tilde F_q$ serves as a fundamental domain of $G$, 
which is the same as $G_q$ introduced in Introduction
and the above-defined covering group of $\pi_q.$
(This appears in Example 1.3 of \cite{CS98}.) Note that the element
$$
V_q:=T_{2q}W_{2q}\inv
=U_{2q}^2
=\begin{pmatrix}\lambda_{2q}^2-1&-\lambda_{2q}\\ \lambda_{2q}&-1\end{pmatrix}
=\begin{pmatrix}\lambda_q+1&-\lambda_{2q}\\ \lambda_{2q}&-1\end{pmatrix}
$$
is elliptic of order $q$
and that $G_q$ is generated by $T_{2q}$ and $V_q$
(see Figure \ref{fig:fq1} and Figure \ref{fig:fq2}).
Here, we note the following property.

\begin{lem}\label{lem:n}
Let $p, q$ be integers with $p,q\ge 2~($where $q$ may be $\infty)$ and
let $K=G_q\cap G_q^{M_p},$ where $M_p\tau=p\tau.$
Then
$|G_q:K|=|G_q^{M_p}:K|.$
\end{lem}

\begin{pf}
If both of the indices are infinite, there is nothing to show.
Assume $n=|G_q:K|<\infty.$
Then the hyperbolic area $\area(Z)$ of the Riemann orbifold $Z=K\backslash\uhp$
is computed as (see \cite[p.~150]{Beardon:disc} or \cite[p.~13]{Katok:fg})
$$
\area(Z)=n\area(G_q\backslash\uhp)=n\area (\tilde F_q)=2n(\pi-\pi/q).
$$
Since $G_q^{M_p}\backslash\uhp$ has the same hyperbolic area
as $G_q\backslash\uhp,$ we have
$$
|G_q^{M_p}:K|
=\frac{\area(K\backslash\uhp)}{\area(G_q^{M_p}\backslash\uhp)}
=\frac{\area(Z)}{\area(G_q\backslash\uhp)}
=n
$$
as required.
When $|G_q^{M_p}:K|<\infty,$ the same argument works for the proof.
\end{pf}

\begin{rem}
In view of the proof, one can also show that the formula
$|G:G\cap G^M|=|G^M:G\cap G^M|$
holds for a cofinite Fuchsian group $G;$ namely, 
if the surface $G\backslash\uhp$ has finite hyperbolic area.
However, this formula is not true in general.
For instance, J\o rgensen, Marden and Pommerenke \cite{JMP81} constructed
a Fuchsian group $\Gamma$ and its subgroup $G$ of index 2 such that
$G^M=M\inv GM$ is a proper subgroup of $G$ for some $M\in\Gamma$
so that $G\cap G^M=G^M.$
In particular, $|G:G\cap G^M|>|G^M:G\cap G^M|=1.$
\end{rem}

The index $|G_q:K|$ can be computed explicitly when $q=3,4,\infty$
(see \cite{Alam21}).
In the present paper, our main concerns are about the limiting case as
$q\to\infty$ and the case of $q=3.$
Those cases will be discussed in Sections \ref{sec:2} and \ref{sec:3}, respectively.

\section{A result on Fuchsian groups}\label{sec:Fuchs}
Let $\Gamma$ be a Fuchsian group acting on the upper half-plane $\uhp$
and $X$ be the quotient Riemann surface $\Gamma\backslash\uhp.$
We denote by $\pi$ the canonical projection $\pi:\uhp\to X.$
When $\gamma\in\Gamma$ is a parabolic element with fixed point at
$\tau_0\in\partial\uhp,$ it is known that $\pi(\tau)$ tends to a puncture,
say $P,$ of $X$ as $\tau\to\tau_0$ nontangentially.
As a convention, we will write $\pi(\tau_0)=P$ in the sequel.
A point $\tau_0\in\uhp$ is called a fixed point of
an element $\gamma$ of $\PSL(2,\R)$ if $\gamma \tau_0=\tau_0$
and  the set of fixed points of $\gamma$ in $\uhp$ is denoted by $\fix(\gamma).$
A non-identity element $\gamma$ is called elliptic if $\fix(\gamma)\ne\emptyset.$
Let $M\in\PSL(2,\R)=\Aut(\uhp)$ such that $M\notin\Gamma.$
Then we consider the number (possibly $\infty$) defined by
\be\label{eq:N}
N_{X_0}(M,\Gamma)=\sup_{x\in X_0}\#\{\pi(M\tau): \tau\in\pi\inv(x)\}.
\ee
for $X_0\subset X.$
The next result is our main lemma and it gives a criterion for finiteness of $N_{X_0}(M,\Gamma).$

\begin{lem}\label{lem:N}
Let $\Gamma$ be a Fuchsian group and $M\in\PSL(2,\R)\setminus\Gamma.$
Suppose that $X_0$ is an uncountable subset of the quotient Riemann surface
$X=\Gamma\backslash\uhp.$
Then $N_{X_0}(M,\Gamma)<\infty$ if and only if $\Gamma^M:=M\inv\Gamma M$
is commensurable with $\Gamma$ in the sense that $|\Gamma:G|<\infty$
for $G=\Gamma\cap \Gamma^M.$
Moreover, in this case, $N_{X_0}(M,\Gamma)=|\Gamma:G|.$
\end{lem}

\begin{pf}
The ``if\," part is almost trivial.
Let $N:=|\Gamma:G|<\infty.$
Then we take $\gamma_1,\dots,\gamma_N\in\Gamma$
so that $\Gamma=G\gamma_1\cup\dots\cup G\gamma_N.$
For $x\in X$ and $\tau_0\in\pi\inv(x),$ we observe
$$
M(\pi\inv(x))=M\Gamma \tau_0
=\bigcup_{j=1}^N MG\gamma_j \tau_0
=\bigcup_{j=1}^N MGM\inv\cdot M\gamma_j \tau_0.
$$
Since $MGM\inv\subset\Gamma,$ the set
$MGM\inv\cdot M\gamma_j \tau_0$ is projected to the point
$\pi(M\gamma_j \tau_0)$ by $\pi.$
Therefore, $\#(\pi(M(\pi\inv(x))))\le N.$
Hence, we conclude that
\be\label{eq:NG}
N_{X_0}(M,\Gamma)
\le |\Gamma:G|.
\ee
To show the ``only if\," part,
we assume that $N:=N_{X_0}(M,\Gamma)$ is finite.
Let
$$
E_1=\bigcup_{\gamma,\delta}\fix(\gamma\inv\delta),
$$
where $\gamma$ and $\delta$ range over
$\gamma\in\Gamma$ and $\delta\in\Gamma^M$ with $\gamma\ne\delta.$
Since each $\fix(\gamma\inv\delta)$ contains at most one point,
the set $E_1$ is at most countable.
Moreover, the set 
$E=\Gamma\cdot E_1=\{\gamma \tau: \gamma\in\Gamma, \tau\in E_1\}$
is also at most countable.
Take a point $\tau_0$ from the uncountable subset $\pi\inv(X_0)\setminus E$ 
of $\uhp$ and fix it.
We regard $G\backslash\Gamma$ as the set of right cosets 
$\{G\gamma: \gamma\in\Gamma\}.$
As we saw above, each set $G\gamma\cdot \tau_0$ projects to one point
$\pi(M\gamma \tau_0)$ under the mapping $\pi\circ M:\uhp\to X.$
We now show that the mapping $\phi:G\backslash\Gamma\to X$
defined by $\phi(G\gamma)=\pi(M\gamma \tau_0)$ is injective.
To this end, we suppose that $\phi(G\gamma_1)=\phi(G\gamma_2)$
for some $\gamma_1,\gamma_2\in\Gamma.$
Then $M\gamma_2 \tau_0=\gamma M\gamma_1 \tau_0$ for some $\gamma\in\Gamma.$
It says that $\gamma_1\tau_0$ is a fixed point of $(\gamma_2\gamma_1\inv)\inv\delta,$
where $\delta=M\inv\gamma M\in\Gamma^M.$
Since $\gamma_1\tau_0\notin E,$ the element $(\gamma_2\gamma_1\inv)\inv\delta$
must be the identity, which implies 
$\gamma_2\gamma_1\inv=\delta\in \Gamma\cap\Gamma^M=G.$
Hence $G\gamma_1=G\gamma_2.$

We have seen that $\phi:G\backslash\Gamma\to X$ is injective.
On the other hand, the image $\phi(G\backslash\Gamma)=\pi(M\Gamma \tau_0)$
consists of at most $N$ points.
Therefore, we obtain
$$
|\Gamma:G|=\# (G\backslash\Gamma)\le N.
$$
Combining with \eqref{eq:NG}, we obtain $|\Gamma:G|=N_{X_0}(M,\Gamma)$
as required.
\end{pf}

In particular, we see that $N_{X_0}(M,\Gamma)$ does not depend on
the uncountable set $X_0\subset X.$ Thus we denote by $N(M,\Gamma)$
the common number $N_{X_0}(M,\Gamma).$

We are now ready to prove Theorem \ref{thm:Fuchs}.

\begin{pf}[Proof of Theorem \ref{thm:Fuchs}]
Recall that $a=(q-1)/(2q)$ and $M_p\tau=p\tau.$
In this proof, we will use the notation introduced in Section \ref{sec:basic}.
By definition, the generalized modular equation \eqref{eq:gme} may be expressed
by $f_a(\beta)=p f_a(\alpha).$
Therefore, 
$\alpha$ and $\beta$ in $\sphere\setminus\{0,1\}$
satisfy \eqref{eq:gme} if and only if $\alpha=\pi_q(\tau), \beta=\pi_q(p\tau)$
for some $\tau\in\uhp.$

First we assume that $\alpha$ and $\beta$ in \eqref{eq:gme} satisfy the
algebraic equation $P(\alpha,\beta)=0,$ where 
$$
P(x,y)=\sum_{j=0}^n a_j(x)y^j,\quad a_j(x)\in \C[x]~(j=0,1,\dots, n),
$$
is an irreducible polynomial of degree $n$ in $\C(x)[y].$
Let $X_0$ be the set of those points $x\in X$ for which $a_k(x)\ne 0$ for some $k.$
Then for a fixed $x_0\in X_0,$ the algebraic equation $P(x_0,y)=0$ in $y$ has at most
$n$ solutions.
Thus we conclude that $N_{X_0}(M_p,G_q)\le n.$
The previous lemma now implies that $|G_q: G_q\cap G_q^{M_p}|\le n.$
By the irreducibility of $P(x,y),$ we see that equality holds.

Conversely, we assume that the equality $|G_q: K|=n$ holds, where
$K=G_q\cap G_q^{M_p}.$
Note that $|G_q^{M_p}: K|=n$ by Lemma \ref{lem:n}.
We denote by $Z$ the quotient Riemann surface $K\backslash\uhp.$
Let $\rho:\uhp\to Z$ be the canonical projection and let $\varphi:Z\to X$
and $\psi:Z\to X$ be the induced mappings satisfying the relations
$\pi_q=\varphi\circ \rho$ and $\pi_q\circ M_p=\psi\circ \rho,$ respectively.
Thus we have the following commutative diagram:
\begin{center}
    \begin{tikzcd}
&\uhp \arrow[r, "M_p"] \arrow[d, "\rho"]\arrow[swap]{dl}{\pi_q}
& \uhp \arrow[d, "\pi_q"] \\
X
&Z \arrow[r, "\psi", labels=below]\arrow[l, "\varphi"]
& X.
\end{tikzcd}
\end{center}

Note that the solution $(\alpha,\beta)$ to the modular equation \eqref{eq:gme}
is now parametrized by $\alpha=\varphi(z)$ and $\beta=\psi(z)$ for $z\in Z.$
We denote by $\hat X$ and $\hat Z$ the compact Riemann surfaces
obtained by filling in the punctures of $X$ and $Z,$ respectively.
Note that $\sphere$ can be taken as $\hat X.$
Then $\varphi$ and $\psi$ extend to $n$-sheeted branched covering maps of $\hat Z$
onto $\hat X=\sphere.$
In particular, $\psi$ may be regarded as a meromorphic function on $\hat Z.$
Let $X_0$ denote $X$ minus the set of critical values of $\varphi:Z\to X;$
indeed $X_0=\sphere\setminus\{0,1,\infty\}$ in this case.
Similarly, $X_1$ is defined for $\psi:Z\to X.$
For a point $x_0\in X_0$ we choose a small disk $U=U(x_0)$ 
with $x_0\in U\subset X_0.$
Let $s_j$ be the elementary symmetric functions of
$\psi\circ \eta_1,\dots, \psi\circ\eta_n$ of degree $j,$ where $\eta_1,\dots,\eta_n$
are local inverses of $\varphi$ of $U.$
Then for each $j=1,2,\dots, n,$ all $s_j:U(x_0)\to X$ piece together to one function and
it extends meromorphically to $\hat X=\sphere$
and $\psi$ satisfies the equation
$$
\psi^n-s_1\circ\varphi\cdot  \psi^{n-1}+\dots 
+(-1)^{n-1}s_{n-1}\circ\varphi\cdot  \psi+(-1)^ns_n\circ\varphi=0
$$
(see \cite[Theorem 8.3]{Forster:RS} for details).
Note that each $s_j:\sphere\to\sphere$ is a rational function.
By writing $(-1)^js_j(x)=a_j(x)/a_0(x)$ for polynomials
$a_0(x),\dots, a_n(x)\in \C[x]$ without non-trivial common factor,
we define 
$$
P(x,y)=a_0(x)y^n+a_1(x)y^{n-1}+\dots+a_{n-1}(x)y+a_n(x).
$$
Then $P(\alpha,\beta)=P(\varphi(z), \psi(z))=0$ for $z\in Z.$
Conversely, suppose that $x_0\in X_0$ with $a_0(x_0)\ne0$
and $y_0\in X_1$ satisfy $P(x_0,y_0)=0.$
Let $\varphi\inv(x_0)=\{z_1,\dots, z_n\}$ and
let $\sigma_1,\dots,\sigma_n$ be the elementary symmetric functions
of $\psi(z_1),\dots,\psi(z_n).$
Note that $\sigma_j=s_j(x_0)$ for $j=1,\dots, n.$
Then,
\begin{align*}
a_0(x_0)\big(y_0-\psi(z_1)\big)\cdots \big(y_0-\psi(z_n)\big)
&=a_0(x_0)\big(y_0^n-\sigma_1 y_0^{n-1}+\cdots
+(-1)^n \sigma_n \big) \\
&=\sum_{j=0}^n a_j(x_0)y_0^{n-j}=P(x_0,y_0)=0.
\end{align*}
Hence, $y_0=\psi(z_k)$ for some $k.$
Since $x_0=\varphi(z_k), y_0=\psi(z_k),$
we have shown that the converse is true.
We next show that the polynomial $P(x,y)$ is irreducible in $y$ with coefficients in $\C[x].$
Suppose, on the contrary, that $P(x,y)$ reduces to the product
$P_1(x,y)P_2(x,y)$ of nonconstant polynomials $P_1, P_2.$
Let $Y_l=\{(x,y)\in X_0\times X_1: P_l(x,y)=0\}$ for $l=1,2.$
We claim that $Y_1\cap Y_2=\emptyset.$
Indeed, if $(x_0, y_0)\in Y_1\cap Y_2,$ then the polynomial
$P(x_0, y)=P_1(x_0,y)P_2(x_0,y)$ has a multiple zero at $y=y_0,$
which is impossible because the set $\psi\inv(y_0)$ consists of $n$
points for $y_0\in X_1.$
Since $\varphi\times\psi:Z_0\to X_0\times X_1$ is continuous,
where $Z_0=\varphi\inv(X_0)\cap\psi\inv(X_1),$ and $Z_0$ is connected,
the image $(\varphi\times\psi)(Z_0)$ is contained in either $Y_1$ or $Y_2.$
However, this is impossible.
The claim has been shown.

Finally, consider the polynomial $P(x,y)$ in $\C(y)[x]$. Then, $P(x,y)$ is a polynomial in $x$ with coefficients in $\C[y],$ and
we may write
$$
P(x,y)=b_0(y)x^m+b_1(y)x^{m-1}+\cdots +b_{m-1}(y)x+b_m(y),
\quad b_j[y]\in\C[y],
$$
where $\displaystyle m=\max_{0\le j\le n} \deg a_j(x).$
Note that $b_0\ne 0.$
Since $\psi:Z\to X$ is $n$-sheeted, as in the case of $\varphi:Z\to X,$
we can see that $m=n$ and $P(x,y)$ is irreducible in $\C(y)[x].$
\end{pf}

Let us summarize the above observations for later use.
Suppose that $K=G_q\cap G_q^{M_p}$ is a subgroup of $G_q$ of finite index $n.$
The intermediate Riemann surface $Z=K\backslash\uhp$ may be used as
a parameter space of the solutions $(\alpha,\beta)$ to \eqref{eq:gme}.
Indeed, the solutions are given by $\alpha=\varphi(z), \beta=\psi(z)$
for $z\in Z,$ where $\varphi, \psi: Z\to X=G_q\backslash\uhp$ 
are (possibly branched) covering
maps satisfying the relations $\varphi(\rho(\tau))=\pi_q(\tau)$
and $\psi(\rho(\tau))=\pi_q(p\tau)$ for $\tau\in\uhp.$
Note that $\varphi$ and $\psi$  extend to the compactifications
$\hat Z$ to $\hat X=\sphere$ as $n$-sheeted branched (analytic) covering maps.
The polynomial $P(x,y)$ whose zero set describes the solutions $(\alpha,\beta)$
can be computed as in the above proof.
The following lemma helps us to compute $\psi$ when we know about $\varphi.$

\begin{lem}\label{lem:omega}
Under the assumption of Theorem \ref{thm:Fuchs}, the M\"obius transformation
$SM_p\tau=-1/(p\tau)$ induces an analytic involution $\omega:Z\to Z$ 
which satisfies the relation  $\rho\circ SM_p=\omega\circ\rho$ and
the functional equation
$$
\psi=1-\varphi\circ\omega.
$$
\end{lem}

\begin{pf}
Let $K=G_q\cap G_q^{M_p}.$
We first note that the M\"obius transformations $S$ and $M_p$ satisfy
the relation $SM_p=M_p\inv S$ and thus $M_pSM_p=S.$
Recall that $S$ normalizes $G_q;$ namely, $G_q^S:=S\inv G_qS=G_q.$
Therefore,
$$
K^{SM_p}=G_q^{SM_p}\cap (G_q^{M_p})^{SM_p}
=G_q^{M_p}\cap G_q^{M_pSM_p}
=G_q^{M_p}\cap G_q^{S}=K,
$$
which means that $SM_p:\uhp\to\uhp$ descends to
an automorphism $\omega:Z\to Z$ such that 
$\rho\circ SM_p=\omega\circ\rho.$
Since $(SM_p)^2=I,$ we see that $\omega\circ\omega=\id.$
We recall Lemma \ref{lem:recip} which says that $\pi_q\circ S=1-\pi_q.$
Then we compute
$$
\psi\circ\omega\circ\rho=\psi\circ\rho\circ SM_p
=(\pi_q\circ M_p)\circ SM_p
=\pi_q\circ S
=1-\pi_q=1-\varphi\circ\rho
$$
and therefore we have $\psi\circ\omega=1-\varphi.$
Since $\omega$ is an involution, we obtain the required relation.
\end{pf}

\begin{rem}
Since $\omega$ comes from the normalizer of $K,$
$\omega$ is indeed an automorphism of the Riemann \emph{orbifold} $K\backslash\uhp.$
In particular, $\omega$ maps a cone point of angle $2\pi/m$ to another 
(possibly the same) cone point of the same angle for an integer $m\ge 2.$
\end{rem}

\section{Applications to the classical case $q=\infty$}\label{sec:2}
As $q\to \infty,$ the Hecke group $H_{2q}$ converges to the group
$H_\infty=\langle S, T_\infty\rangle,$ which is a subgroup of
the modular group $\PSL(2,\Z)$ of index 3, where $T_\infty z=z+2.$
At the same time, the subgroup $G_q$ converges to
the principal congruence subgroup $\Gamma(2)$ of level 2, which will
be denoted by $G_\infty$ in this section.
That is to say,
$$
G_\infty=\{A\in \SL(2,\Z): A\equiv I \mod 2\}/\{\pm I\}.
$$
This corresponds to the case when $a=1/2$ in the modular equation \eqref{eq:gme}. The fundamental domain $\tilde F_\infty$ is described by $0\le\Re \tau\le 2,
|\tau-1/2|\ge 1/2, |\tau-3/2|\ge 1/2$ for $\tau\in\uhp.$
A pair of generators of $G_\infty$ are given by
$$
T=T_\infty=\begin{pmatrix}1&2\\0&1\end{pmatrix}
\aand
V=\begin{pmatrix}3&-2\\ 2&-1\end{pmatrix}.
$$
Note that $V$ is a parabolic element with fixed point at $\tau=1$ and that
$V(0)=2.$
Also, let
$$
W=-V\inv T=\begin{pmatrix}1&0\\2&1\end{pmatrix}.
$$
It is well known that $\Gamma(2)=G_\infty$ is a Fuchsian group
uniformizing the thrice punctured sphere $\sphere\setminus\{0,1,\infty\}$
(see \S\S 3.4-5 in Chapter 7 of \cite{Ahlfors:ca} and \S\S 4.3 in Chapter 1 of \cite{Farkas:RS}, where the symbol $\lambda$
is used for $\pi_\infty$).
The canonical projection $\pi_\infty:\uhp\to X$ may be considered to be
$\pi_\infty:\uhp\to \sphere\setminus\{0,1,\infty\}=\C\setminus\{0,1\}$
with 
\be\label{eq:pi}
\pi_\infty(0)=1, \quad \pi_\infty(\infty)=0 \aand
\pi_\infty(1)=\infty
\ee
by \eqref{eq:fa}.
Let $\rho:\uhp\to Z, \varphi,\psi:Z\to X$ be as in the proof of Theorem \ref{thm:Fuchs}.
Here is a simple but useful observation.
The conjugation of a matrix by the M\"obius transformation $M_p$ is computed by
$$
M_p\inv \begin{pmatrix}a& b\\ c&d\end{pmatrix}M_p
=\begin{pmatrix}a& b/p\\ pc&d\end{pmatrix}.
$$
Hence, an element $A=\begin{pmatrix}a& b\\ c&d\end{pmatrix}$
of $G_\infty$ is a member of $G_\infty^{M_p}$ precisely if
$c\equiv 0\,(\mod 2p).$

\subsection{Case $p=2$}
\begin{figure}[h]
         \centering
         \includegraphics[width=0.7\textwidth]{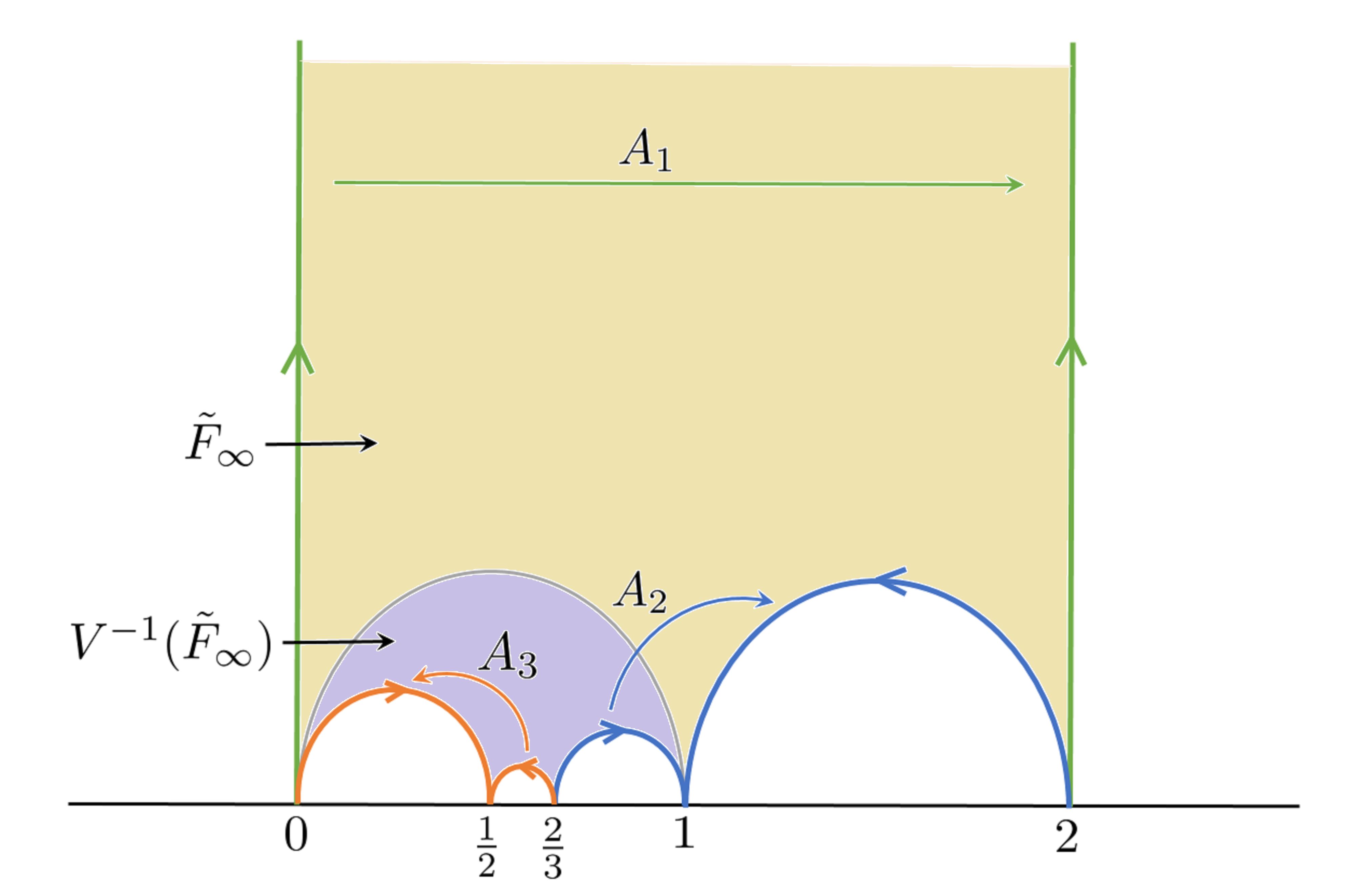}
         \caption{Fundamental domain $F$ of $K=G_\infty\cap G_\infty^{M_2}$.}
         \label{fig:fd1}
\end{figure}

We first look at the most classical case when $p=2.$
Let $K=G_\infty\cap G_\infty^{M_2}.$
Then 
$$
A_1:=T=\begin{pmatrix}1&2\\0&1\end{pmatrix}, \quad
A_2:=V^2=\begin{pmatrix}5& -4\\ 4&-3\end{pmatrix}, \quad
A_3:=V\inv TV=\begin{pmatrix}-3& 2\\ -8&5\end{pmatrix}
$$
are side-paring transformations of the hyperbolic polygon
$F=\tilde F_\infty\cup V\inv(\tilde F_\infty)$
({see Figure \ref{fig:fd1}}).

Therefore, the elements $A_1, A_2, A_3$ generate a subgroup
of $G_\infty$ of index 2.
In view of the forms of $A_j,$ we see that $A_j\in K$ for $j=1,2,3$
and thus we have $K=\langle A_1, A_2, A_3\rangle$ and
$|G_\infty:K|=2.$
It is easy to observe that the quotient surface $Z=K\backslash\uhp$
is conformally equivalent to a four-times punctured sphere.
Thus we can assume that $Z=\sphere\setminus\{0, \infty, 1, b \},$
where $\rho(0)=0, \rho(1)=\infty, \rho(\infty)=1, \rho(1/2)=b.$
In view of \eqref{eq:pi}, we have the conditions
$\varphi(0)=1, \varphi(\infty)=\infty, \varphi(1)=\varphi(b)=0.$
Here we used the fact that $1/2$ is conjugate to $\infty$ under the
action of $G_\infty.$
We note that the extension
$\varphi:\hat Z=\sphere\to \hat X=\sphere$ is a rational map
of degree 2.
Since $\varphi\inv(1)=\{0\}, \varphi\inv(\infty)=\{\infty\},$
we see that $\varphi$ takes the values $0$ and $\infty$ with multiplicity 2.
Therefore, $\varphi$ should have the form $\varphi(z)=cz^2+1$
for a constant $c\ne0.$
Since $\varphi(1)=0,$ we conclude that $c=-1$ and $b=-1.$
Hence, $\varphi(z)=1-z^2$ under the above normalization.

Next we determine the form of $\psi:Z\to X.$
By the relation $\pi_\infty(2\tau)=\psi(\rho(\tau)),$ we have
the necessary conditions $\psi(0)=\pi_\infty(0)=1, 
\psi(\infty)=\pi_\infty(2)=\pi_\infty(0)=1, \psi(1)=\pi_\infty(\infty)
=0$ and $\psi(-1)=\pi_\infty(1)=\infty.$
In particular, $\psi\inv(0)=\{1\}, \psi\inv(\infty)=\{-1\}.$
Since $\psi:\hat Z=\sphere\to \hat X=\sphere$ is a rational map of degree 2,
$\psi$ has the form
$$
\psi(z)=\frac{c(z-1)^2}{(z+1)^2}
$$
for a constant $c\ne0.$
Since $\psi(0)=1,$ we have $c=1.$
We note that $\psi(\infty)=1$ is also satisfied.
In this way, the solution $(\alpha,\beta)$ of the modular equation
\eqref{eq:gme} with $a=1/2, p=2$ is parametrized by
$$
\alpha=\varphi(z)=1-z^2
\aand
\beta=\psi(z)=\frac{(z-1)^2}{(z+1)^2}.
$$
Eliminating the variable $z,$ we see that
$$
\beta=\left(\frac{1-\sqrt{1-\alpha}}{1+\sqrt{1-\alpha}}\right)^2,
$$
which is equivalent to \eqref{eq:classical} with $\alpha=r^2$ and $\beta=s^2.$
The polynomial in Theorem \ref{thm:Fuchs} is given by
$$
P(x,y)=x^2y^2-2(x^2-8x+8)y+x^2.
$$
In this case, the polynomial is not symmetric in $x$ and $y.$

\subsection{Case $p=3$}
\begin{figure}[h]
         \centering
         \includegraphics[width=0.8\textwidth]{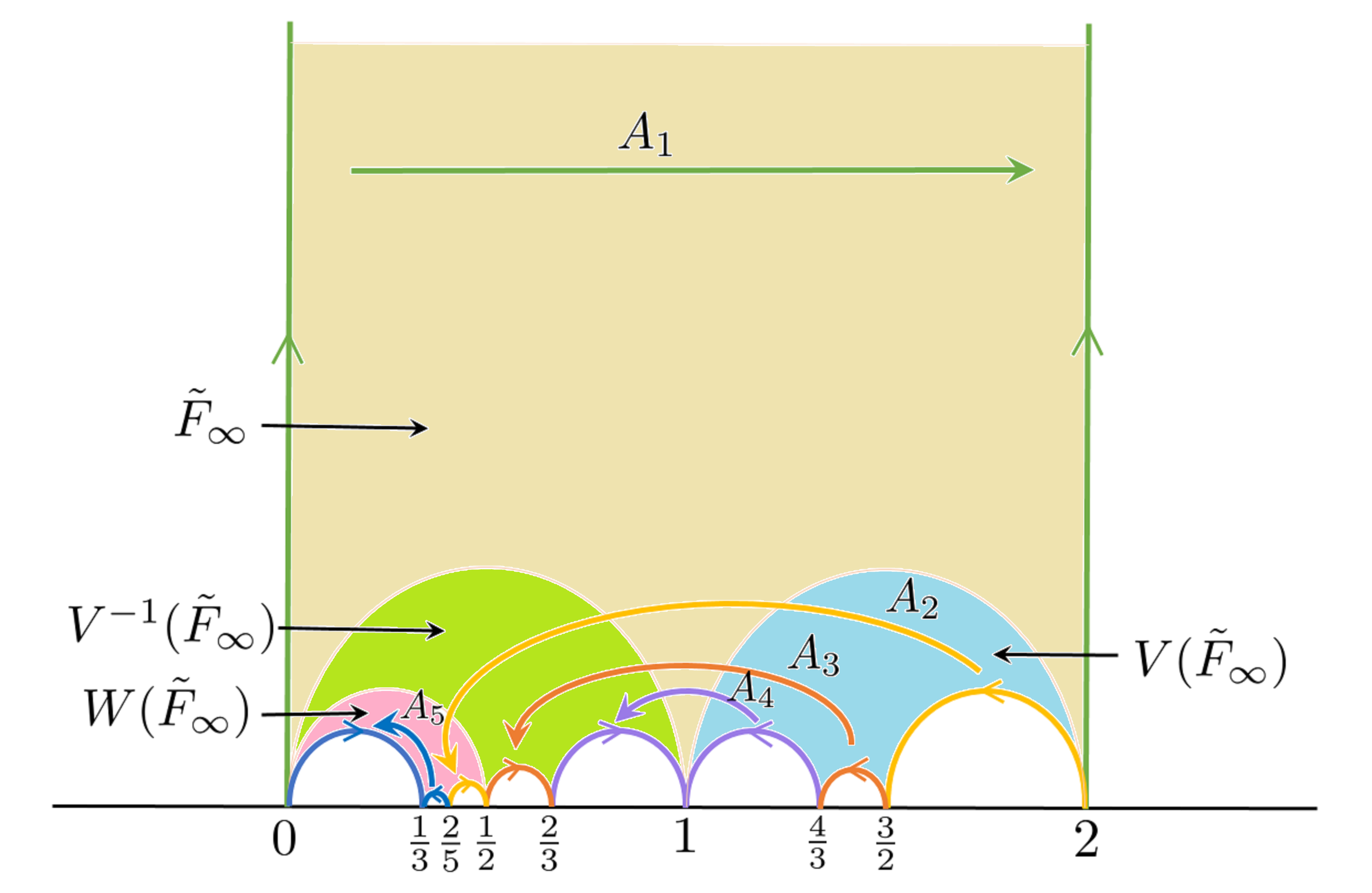}
         \caption{Fundamental domain $F$ of $K=G_\infty\cap G_\infty^{M_3}$.}
         \label{fig:23}
\end{figure}

Next we consider the case when $p=3.$
We now consider the hyperbolic 10-gon
$F=\tilde F_\infty\cup V(\tilde F_\infty)\cup 
V\inv(\tilde F_\infty)\cup W(\tilde F_\infty)$
with vertices at $0, 1/3, 2/5, 1/2, 2/3, 1, 4/3,$ $3/2, 2, \infty$ 
in the counterclockwise order.
Then the elements
\begin{align*}
A_1&:=T,\quad
A_2:=V\inv T^2V\inv=\begin{pmatrix}5&-8\\12&-19\end{pmatrix}, \quad
A_3:=V\inv T\inv V\inv=\begin{pmatrix}-7&10\\-12&17\end{pmatrix}, \quad \\
A_4&:=V^{-3}=\begin{pmatrix}-5&6\\-6&7\end{pmatrix}, \quad
A_5:=V\inv T V\inv T\inv V=\begin{pmatrix}-5&2\\-18&7\end{pmatrix}
\end{align*}
form side-pairing transformations of $F$ ({see Figure \ref{fig:23}}).

Therefore, $K=\langle A_1, A_2, A_3, A_4, A_5\rangle$ is a subgroup of $G_\infty$
for which $F$ is a fundamental domain.
In view of the forms of $A_j,$ we observe that 
$K\subset G_\infty\cap G_\infty^{M_3}.$
On the other hand, the elements $V, V\inv, W$ are not contained in 
$G_\infty^{M_3},$ therefore $K=G_\infty\cap G_\infty^{M_3}$
and $|G_\infty:K|=4.$
It is a simple task to see that the quotient surface $Z=K\backslash\uhp$
is a six-times punctured sphere.
Recall that $\varphi,\psi:Z\to X$ extend to rational functions of degree 4.
We may normalize the punctures so that
$\rho(\infty)=0, \rho(0)=\rho(2/5)=\rho(2)=1, \rho(1/3)=\infty,
\rho(1/2)=\rho(3/2)=b, \rho(1)=c, \rho(2/3)=\rho(4/3)=d.$
Thus $\varphi(0)=\pi_\infty(\infty)=0.$
Similarly, we obtain $\varphi(1)=1, \varphi(\infty)=\infty, 
\varphi(b)=0, \varphi(c)=\infty, \varphi(d)=1.$
We have to compute multiplicities of these values.
For instance, we have $\varphi\inv(\infty)=\{\infty, c\}.$
The part of the basic fundamental domain $\tilde F_\infty$ corresponding to
$\infty$ under $\pi_\infty$ is the cusp neighbourhood 
$D=\tilde F_\infty\cap\{\tau\in\uhp: |\tau-1|<\ve\}$
for a small enough $\ve>0.$
Since $V$ and $V\inv$ fix $1$ while $W$ sends $1$ to $1/3,$
the multiplicity of $\varphi$ at $\rho(1)=c$ is 3 and that is 1 at $\rho(1/3)=\infty.$
We write $(\varphi)_\infty=3\cdot c+1\cdot\infty$ as a divisor for short\footnote{%
When the equation $\varphi(z)=w$ has solutions $z_j$ with multiplicities $m_j$
for $j=1,2,\dots, N,$ we write $(\varphi)_w=m_1\cdot z_1+m_2\cdot z_2+\cdots
+m_N\cdot z_N$ as a divisor on $Z$.}.
In the same way, we have $(\varphi)_0=3\cdot b+1\cdot 0$
and $(\varphi)_1=1\cdot d+3\cdot 1.$
In particular, $\varphi$ may be expressed by
$$
\varphi(z)=e\frac{z(z-b)^3}{(z-c)^3}
\aand
\varphi(z)-1=e'\frac{(z-d)(z-1)^3}{(z-c)^3}
$$
for constants $e$ and $e'.$
Since $b,c,d$ are different from $1,$ we have the unique solution
$b=-2, c=-1/2, d=-1$ and $e=e'=1/8.$
Hence,
$$
\alpha=\varphi(z)=\frac{z(z+2)^3}{(2z+1)^3}.
$$
Next we determine the form of $\omega$ given in Lemma \ref{lem:omega}.
Since $SM_3$ swaps $0$ and $\infty$ (respectively $1/3$ and $ -1$),
$\omega$ swaps $\rho(0)=1$ and $\rho(\infty)=0$
(respectively, $\rho(1/3)=\infty$ and $\rho(-1)=\rho(1)=-1/2$).
Hence, $\omega(z)=(1-z)/(1+2z).$
By Lemma \ref{lem:omega}, we obtain
$$
\beta=\psi(z)=1-\varphi(\omega(z))=\frac{z^3(z+2)}{2z+1}.
$$
We compute
$$
\alpha\beta=\frac{z^4(z+2)^4}{(2z+1)^4}
\aand
(1-\alpha)(1-\beta)=\frac{(1-z^2)^4}{(2z+1)^4}.
$$
Note that $\varphi$ and $\psi$ both map the interval $[0,1]$
onto itself homeomorphically.
Hence, for $\alpha,\beta\in[0,1],$ we obtain the relation

\begin{equation}\label{eq:le}
    (\alpha\beta)^{1/4}+\big\{(1-\alpha)(1-\beta)\big\}^{1/4}=1.
\end{equation}

This is known as Legendre's modular equation
(see (5.37) in \cite{AVV:conf} and (4.1.16) in \cite{BB:AGM}). 
Formula \eqref{eq:le} may be transformed to the polynomial equation 
$P(\alpha, \beta)=0$ as in Theorem \ref{thm:Fuchs}, where
$$
P(x,y)=y^4+2x^3y^3-2xy-x^4.
$$
This is known as the modular equation of third order in Jacobi's form 
(see (3.42) in \cite{Ranjan}).

\section{Case of degree $3$ in the theory of signature $3$}\label{sec:3}
When $q=3,$ the standard generators of $G_3$ are
$$
T:=T_{6}=\begin{pmatrix}1&\sqrt 3\\0&1\end{pmatrix}
\aand
V:=V_3=\begin{pmatrix}2&-\sqrt 3\\ \sqrt 3&-1\end{pmatrix}.
$$
Recall that $V$ is an elliptic element of order 3 with fixed point at
$\tau_0:=e^{\pi i/6}=(\sqrt 3+i)/2.$
Recall also that the canonical projection 
$\pi_3:\uhp\to X=\sphere\setminus\{0,1\}$ satisfies
\be\label{eq:pi3}
\pi_3(0)=1, \quad \pi_3(\infty)=0 \aand
\pi_3(\tau_0)=\infty
\ee
by \eqref{eq:fa}.
The following result is useful.

\begin{lem}\label{lem:3}
$A\in\PSL(2,\R)$ belongs to $G_3$ precisely when $A$ is represented by a
matrix of the form
\begin{equation}\label{eq:even}
\begin{pmatrix}a & b\sqrt3 \\ c\sqrt 3& d\end{pmatrix}, \quad
a,b,c,d\in\Z,~ad-3bc=1.
\end{equation}
\end{lem}

\begin{pf}
We denote by $G$ the M\"obius transformations represented by
the matrices in \eqref{eq:even}.
Then it is well known that $G$ is a proper subgroup of $H_6$
(see \cite{Hut02}).
On the other hand, it is obvious that $G_3$ is contained in $G.$
Since $|H_6:G_3|=2,$ we have $G=G_3$ as required.
\end{pf}

\subsection{Case $p=3$}
\begin{figure}[h]
         \centering
         \includegraphics[width=0.7\textwidth]{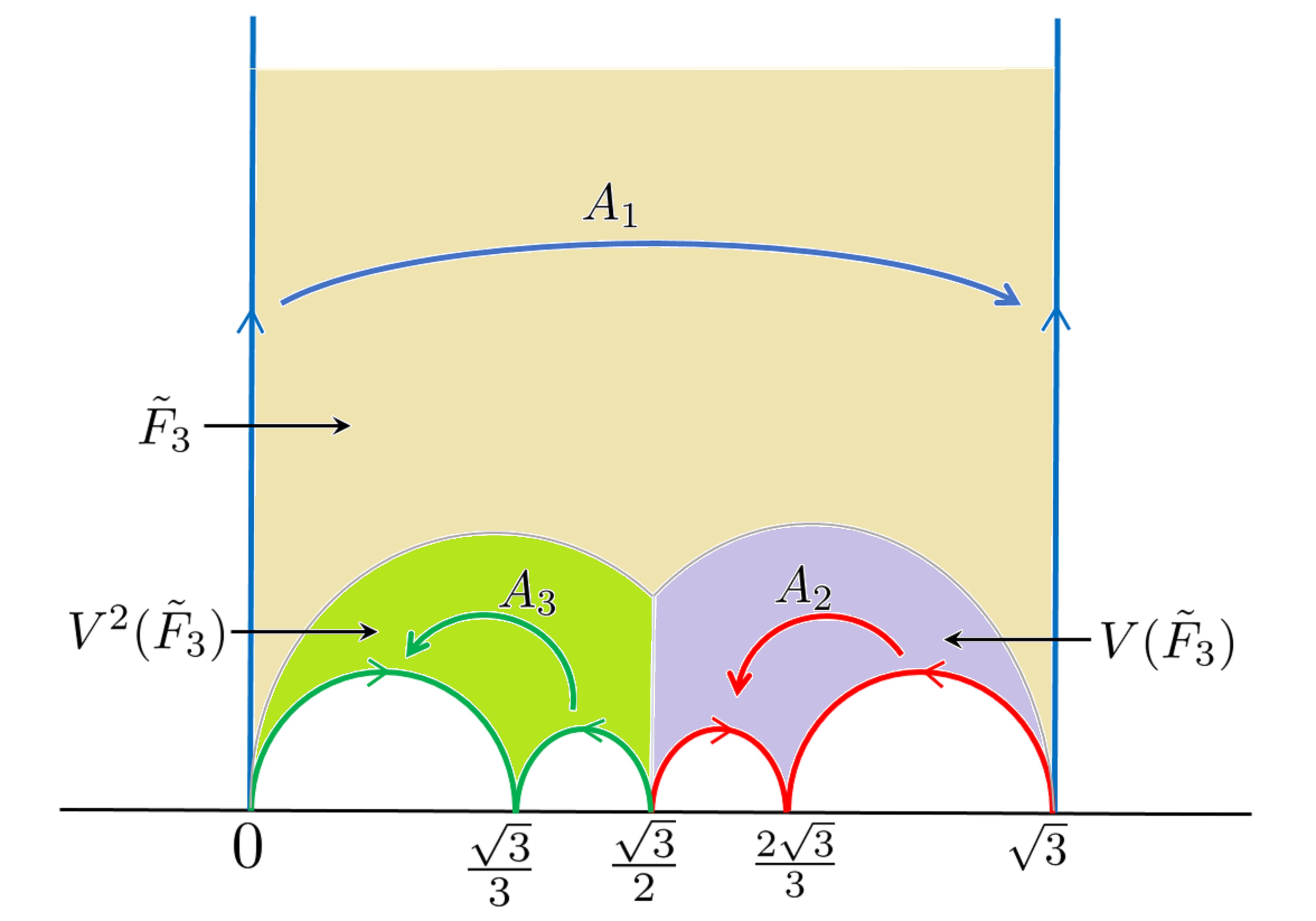}
         \caption{Fundamental domain $F$ of $K=G_3\cap G_3^{M_3}$.}
         \label{fig:33}
\end{figure}

Since the case of $p=3$ is simpler than that of $p=2,$ we start with this case.
Let $F=\tilde F_3\cup V(\tilde F_3)\cup V^2(\tilde F_3).$
Then $F$ is a hyperbolic polygon with six vertices at
$0, \sqrt 3/3, \sqrt 3/2, 2\sqrt 3/3, \sqrt 3, \infty.$
Let $K$ be the group generated by the side pairing transformations
$$
A_1=T, \quad
A_2=VTV\inv=\begin{pmatrix}-5 & 4\sqrt3 \\ -3\sqrt 3& 7\end{pmatrix}, \quad
A_3=V\inv TV=\begin{pmatrix}-2 & \sqrt3 \\ -3\sqrt 3& 4\end{pmatrix}
$$
of $F$ (see Figure \ref{fig:33}).

It is easy to see that $K=G_3\cap G_3^{M_3}$ by Lemma \ref{lem:3} and
that $K$ is torsion-free and satisfies $|G_3:K|=3.$
Since the quotient surface $Z=K\backslash\uhp$ is a four-times punctured sphere,
we may normalize $\rho:\uhp\to Z$ so that 
$\rho(\infty)=1, \rho(\sqrt 3/3)=a, \rho(2\sqrt 3/3)=b,
\rho(0)=\rho(\sqrt 3/2)=\rho(\sqrt 3)=0$ and $\rho(\tau_0)=\infty.$
By the form of $K,$ the element $V$ normalizes $K;$ that is, $V\inv KV=K.$
Thus $V$ induces an analytic automorphism $v:Z\to Z.$
We recall that $V$ is a rotation about $\tau_0$ of angle $-2\pi/3.$
It is clear that $v(0)=0$ and $v(\infty)=\infty.$
Since $V$ satisfies $V(\infty)=2\sqrt 3/3, V(2\sqrt 3/3)=\sqrt 3/3$
and $V(\sqrt 3/3)=\infty,$ the map $v$ should have the form
$v(z)=e^{-2\pi i/3} z.$
Hence, $b=e^{-2\pi i/3}=-(1+i\sqrt 3)/2$ and $a=\bar b=(-1+i\sqrt 3)/2.$
We will determine the forms of rational maps 
$\varphi, \psi:\hat Z=\sphere\to\hat X=\sphere$ of degree 3.
Since $\varphi$ has a branch point at $\infty$ of order 3,
we have $\varphi\inv(\infty)=\{\infty\}.$
In particular, $\varphi$ is a polynomial of degree 3.
We also have $\varphi\inv(1)=\{0\}.$
Therefore, $\varphi$ should be of the form $1+c z^3$
for a constant $c.$
Since $\varphi(1)=\pi_3(\infty)=0,$ we obtain $c=-1.$
Next we determine $\omega$ in Lemma \ref{lem:omega}.
Since $SM_3$ swaps $0$ and $\infty$ (respectively $\sqrt 3/3$
and $-\sqrt 3/3\equiv 2\sqrt 3/3~ (\mod K),$
$\omega$ swaps $1$ and $0$ (respectively, $a$ and $b$),
Hence, after some computations, we get the form $\omega(z)=(1-z)/(1+2z).$
We now apply Lemma \ref{lem:omega} to obtain 
$\psi(z)=1-\varphi(\omega(z))=\omega(z)^3=(1-z)^3/(1+2z)^3.$
In summary, we have
$$
\alpha=\varphi(z)=1-z^3
\aand
\beta=\psi(z)=\frac{(1-z)^3}{(1+2z)^3}.
$$
We are now able to show Theorem \ref{thm:33} easily.
The polynomial $P(x,y)$ in Theorem \ref{thm:Fuchs} is computed as
\begin{align*}
P(x,y)&=(8 x-9)^3 y^3+3 \left(64 x^3+504 x^2-1053 x+486\right) y^2 \\
&\quad +3 \left(8 x^3-171 x^2+405 x-243\right) y+x^3.
\end{align*}

\subsection{Case $p=2$}
\begin{figure}[h]
         \centering
         \includegraphics[width=0.7\textwidth]{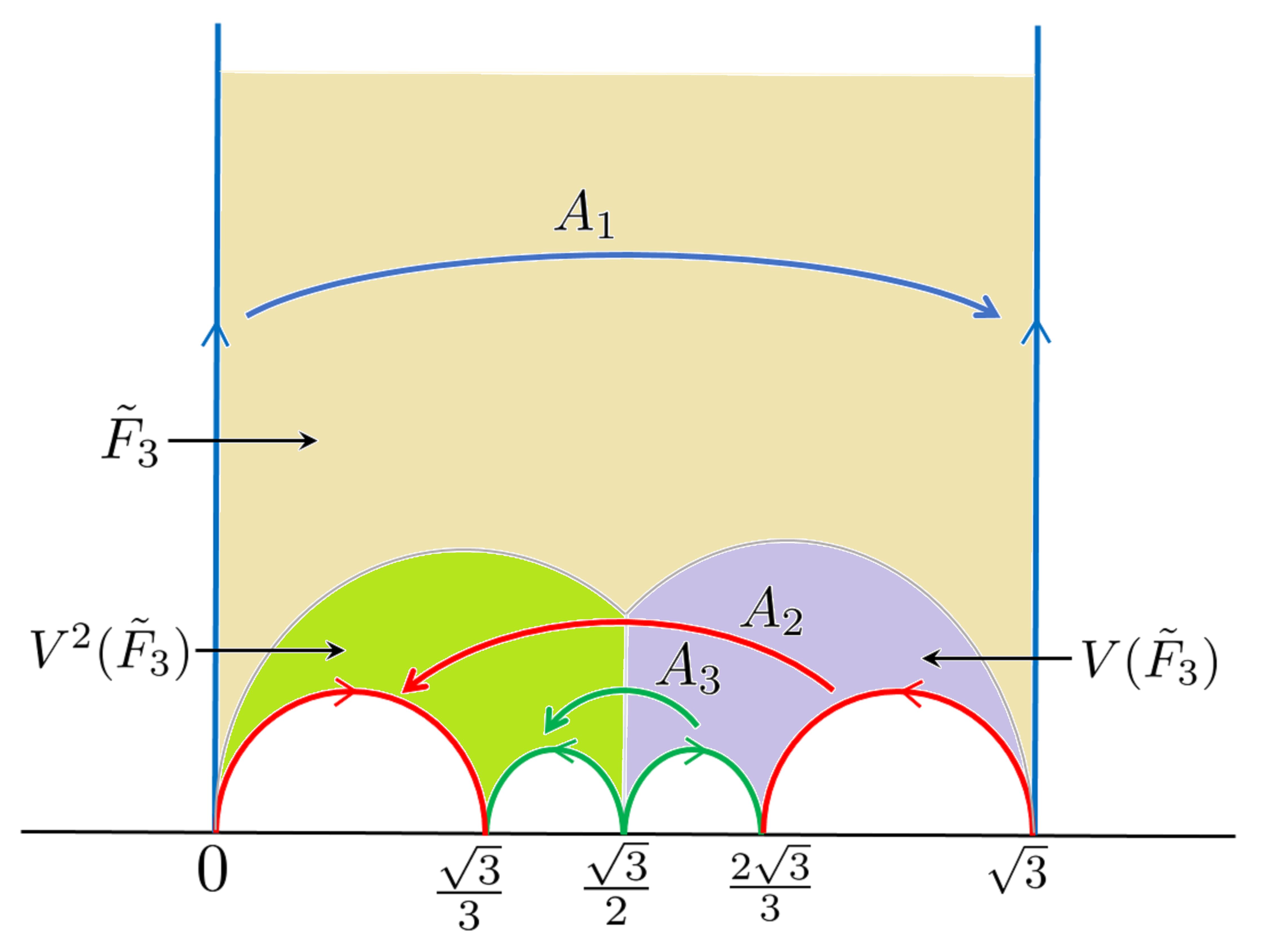}
         \caption{Fundamental domain $F$ of $K=G_3\cap G_3^{M_2}$.}
         \label{fig:32}
\end{figure}

Let $F$ be the same as in the previous subsection.
We choose side-pairing transformations $A_1, A_2, A_3$ of $F$
at this time as follows:
$$
A_1=T,\quad
A_2=V^2TV^2=\begin{pmatrix}1 & -\sqrt3 \\ 2\sqrt 3& -5\end{pmatrix}, \quad
A_3=(SV)\inv TSV=\begin{pmatrix}5 & -3\sqrt3 \\ 4\sqrt 3& -7\end{pmatrix}.
$$
Then $A_1, A_2$ and $A_3$ generate the torsion-free group 
$K=G_3\cap G_3^{M_2}$ and the relation $|G_3:K|=3$ follows (see Figure \ref{fig:32}).

We also see that $Z=K\backslash\uhp$ is a four-times punctured sphere.
We normalize the map $\rho:\uhp\to Z$ so that $\rho(\infty)=0,
\rho(0)=\rho(\sqrt 3)=1$ and $\rho(\tau_0)=-1.$
Let $b=\rho(\sqrt 3/3)=\rho(2\sqrt 3/3)$ and $c=\rho(\sqrt 3/2).$
Since $\pi_3(\tau)=\varphi(\rho(\tau))$ has a branch point of order 3
at $\tau=\tau_0,$ the map $\varphi(z)$ has a branch point of order 3
at $z=-1;$ that is, $\varphi\inv(\infty)=\{-1\}.$
Since $V(\infty)=2\sqrt 3/3, V^2(\infty)=\sqrt 3/3,$
we have $(\varphi)_0=1\cdot 0+2\cdot b.$
Therefore, the rational map $\varphi$ of degree $3$ should have the form
$$
\varphi(z)=\frac{mz(z-b)^2}{(z+1)^3}
$$
for a constant $m\ne0.$
Observe that $V(\tilde F_3)$ and $V^2(\tilde F_3)$ share 
the cusps at $\sqrt 3$ and $0,$ respectively, with $\tilde F_3.$
Hence, we see that $(\varphi)_1=2\cdot 1+1\cdot c$ and
$$
\varphi(z)-1=\frac{m'(z-c)(z-1)^2}{(z+1)^3}
$$
for a constant $m'\ne0.$
Then we have $b=-3, c=-2, m=1/2$ and $m'=-1/2$ so that
$\varphi(z)=z(z+3)^2/[2(z+1)^3].$
Next we will determine $\omega$ in Lemma \ref{lem:omega}.
Since $SM_2$ swaps $0$ and $\infty$ (respectively, $\sqrt 3/3$
and $-\sqrt 3/2\equiv \sqrt 3/2 \mod K$),
$\omega$ swaps $1$ and $0$ (respectively, $b=-3$ and $c=-2$),
the form of $\omega$ is $(1-z)/(1+z).$
Now Lemma \ref{lem:omega} yields $\psi(z)=1-\varphi(\omega(z))
=z^2(z+3)/4.$
We summarize the results as
$$
\alpha=\varphi(z)=\frac{z(z+3)^2}{2(z+1)^3}
\aand
\beta=\psi(z)=\frac{z^2(z+3)}{4}.
$$
We remark that these expressions of $\alpha$ and $\beta$
are found in Ramanujan's notebook (see \cite[Theorem 6.1]{BBG95}).
Since
$$
\alpha\beta=\frac{z^3(z+3)^3}{8(z+1)^3}
\aand
(1-\alpha)(1-\beta)=\frac{(1-z)^3(z+2)^3}{8(z+1)^3},
$$
it is now easy to obtain the formula in Theorem \ref{thm:23}.
We also obtain the polynomial $P(x,y)$ in Theorem \ref{thm:Fuchs} as
$$
P(x,y)=(2x-1)^3y^3-3x(4x^2-13x+10)y^2+3x(2x^2-10x+9)y-x^3,
$$
which is equivalent to \eqref{eq:me23}.


\def\cprime{$'$} \def\cprime{$'$} \def\cprime{$'$}
\providecommand{\bysame}{\leavevmode\hbox to3em{\hrulefill}\thinspace}
\providecommand{\MR}{\relax\ifhmode\unskip\space\fi MR }
\providecommand{\MRhref}[2]{%
  \href{http://www.ams.org/mathscinet-getitem?mr=#1}{#2}
}
\providecommand{\href}[2]{#2}


\begin{thebibliography}{10}

\bibitem{Ahlfors:ca}
L.~V. Ahlfors, \emph{Complex {A}nalysis, 3rd ed.}, McGraw Hill, New York, 1979.

\bibitem{Alam21}
Md.~S. Alam, \emph{On {R}amanujan's modular equations and {H}ecke groups},
  Preprint (2021).

\bibitem{AQVV00}
G.~D. Anderson, S.-L. Qiu, M.~K. Vamanamurthy, and M.~Vuorinen,
  \emph{Generalized elliptic integrals and modular equations}, Pacific J. Math.
  \textbf{192} (2000), 1--37.

\bibitem{ASVV10}
G.~D. Anderson, T.~Sugawa, M.~K. Vamanamurthy, and M.~Vuorinen,
  \emph{Hyperbolic metric on the twice punctured sphere with one conic
  singularity}, Math. Z. \textbf{266} (2010), 181--191.

\bibitem{AVV:conf}
G.~D. Anderson, M.~K. Vamanamurthy, and M.~K. Vuorinen, \emph{Conformal
  {I}nvariants, {I}nequalities, and {Q}uasiconformal {M}aps},
  Wiley-Interscience, 1997.

\bibitem{Beardon:disc}
A.~F.~Beardon, \emph{The Geometry of Discrete Groups},
 Graduate Texts in Mathematics \textbf{91}, Springer-Verlag, New York, 1983.

\bibitem{Berndt: NotebookI}
  B.~C. Berndt, \emph{Ramanujan's Notebooks, Part \RomanNumeralCaps{1}}, Springer-Verlag, New York,
  1985.
  
\bibitem{Berndt: NotebookII}
  B.~C. Berndt, \emph{Ramanujan's Notebooks, Part \RomanNumeralCaps{2}}, Springer-Verlag, New York,
  1989.

\bibitem{Berndt: NotebookIII}
  B.~C. Berndt, \emph{Ramanujan's Notebooks, Part \RomanNumeralCaps{3}}, Springer-Verlag, New York,
  1991.
  
\bibitem{Berndt: NotebookV}
  B.~C. Berndt, \emph{Ramanujan's Notebooks, Part \RomanNumeralCaps{5}}, Springer-Verlag, New York,
  1998.

\bibitem{BBG95}
B.~C. Berndt, S.~Bhargava, and F.~G. Garvan, \emph{Ramanujan's theories of
  elliptic functions to alternative bases}, Trans. Amer. Math. Soc.
  \textbf{347} (1995), 4163--4244.

\bibitem{BB:AGM}
J.~M. Borwein and P.~B. Borwein, \emph{Pi and the {AGM}}, Wiley, New York,
  1987.

\bibitem{Carat:f2}
C. Carath\'eodory, \emph{Funktionentheorie. {B}and {II}},
Birkh\"auser Verlag, Basel, 1950.

\bibitem{CS98}
\.{I}.~N. Cang\"ul and D.~Singerman, \emph{Normal subgroups of {H}ecke groups
  and regular maps}, Math. Proc. Camb. Phil. Soc. \textbf{123} (1998), 59--74.

\bibitem{Farkas:RS}
H.~M. Farkas and I.~Kra, \emph{Theta {C}onstants, {R}iemann {S}urfaces and the {M}odular {G}roup}, Graduate Studies in Mathematics, vol. 37, American Mathematical Society, 2001.

\bibitem{Forster:RS}
O.~Forster, \emph{Lectures on {R}iemann {S}urfaces}, Springer-Verlag, New York,
  1981.

\bibitem{Hut02}
J.~I. Hutchinson, \emph{On a class of automorphic functions}, Trans. Amer.
  Math. Soc. \textbf{3} (1902), 1--11.
  
\bibitem{JMP81}
T.~J\o rgensen, A.~Marden, and C.~Pommerenke, \emph{Two examples of covering
  surfaces}, Riemann surfaces and related topics: {P}roceedings of the 1978
  {S}tony {B}rook {C}onference ({S}tate {U}niv. {N}ew {Y}ork, {S}tony {B}rook,
  {N}.{Y}., 1978), Ann. of Math. Stud., vol.~97, Princeton Univ. Press,
  Princeton, N.J., 1981, pp.~305--317.
  
\bibitem{JG:dessin}
G.~A. Jones and J.~Wolfart, \emph{Dessins d'{E}nfants on {R}iemann {S}urfaces},
  Springer, New York, 2016.
 
\bibitem{Katok:fg}
S.~Katok, \emph{Fuchsian {G}roups}, The University of Chicago Press, Chicago
  and London, 1992.

\bibitem{LV:qc}
O.~Lehto and K.~I. Virtanen, \emph{Quasiconformal {M}appings in the {P}lane,
  2nd {E}d.}, Springer-Verlag, 1973.

\bibitem{Ram: nb}
S.~Ramanujan, \emph{Notebooks (2 volumes)}, Tata Institute of Fundamental Research, Bombay, 1957.

\bibitem{Ram: lnb}
S.~Ramanujan, \emph{The lost notebook and other unpublished papers}, Narosa, New Delhi, 1988.

\bibitem{Ranjan}
R.~Roy, \emph{Elliptic and Modular Functions from Gauss to Dedekind to Hecke}, Cambridge University Press, 2017.

\end{thebibliography}
\end{document}